\documentclass[twocolumn,oneside,a4paper]{article}
\usepackage{latexsym}
\usepackage{graphicx}
\begin{document}
\title{\textbf{The fractional - order controllers:
       Methods for their synthesis and application}}
\author{Ivo Petras  \thanks{
     Department of Informatics and Process Control,  
     BERG Faculty, Technical University of Kosice, B. Nemcovej 3, 
     042 00, Kosice, Slovak  Republic, e-mail:\, \textit{petras@tuke.sk}}}         
\date{}

\maketitle
\noindent
\begin{abstract}
This paper deals with fractional-order  controllers. We outline
mathematical description of fractional controllers and methods of their
synthesis and application. Synthesis method is a modified  root
locus method for  fractional-order systems and fractional-order controllers.
In the next section we describe how to apply the fractional controller
on control systems.
\end{abstract}

\section{Introduction}

The real objects are generally fractional \cite{Torvik}, however, for
many of them the fractionality is very low.  A~typical example of a
non-integer (fractional) order system is the voltage-current relation
of a semi-infinite lossy $RC$ line
or diffusion of the heat into a semi-infinite solid, where heat flow is
equal to the half-derivative of temperature \cite{Podlubny4}:
$$       \frac{d^{0.5}T(t)} {dt^{0.5}}    = q(t). $$
Their integer-order description can cause significant differences
between mathematical model and the real system. The main reason
for using integer-order models was the absence of solution methods for
fractional-order differential equations.
Recently, important achievements
\cite{Axtell, Dorcak2, Matignon, Oldham, Podlubny1, Podlubny3, Woon, Zavada}
were obtained, which enable to take into account real order of dynamic systems.

$PID$ controllers belong to dominating industrial controllers and
therefore are objects of steady effort for improvements of their
quality and robustness. One of the possibilities to improve $PID$
controllers is to use fractional-order controllers with non-integer
derivation and integration parts. For fractional-order systems
the fractional controller $CRONE$ \cite{Oustaloup} has been developed,
$PD^{\delta}$ controller \cite{Dorcak2} and the
$PI^{\lambda}D^{\delta}$ controller \cite{Petras3, Podlubny2} has
been suggested.

Before the actual design of controllers for dynamical systems it is
necessary to identify these systems \cite{Dorcak1, Podlubny3}, and then
determine their dynamical properties. The dynamical properties of the system
under observation can be expressed mathematically and graphically in the form
of various characteristics. These characteristics can be determined from the
differential equation or transfer function via computation or experimentally
by exciting the system from the equilibrium.
Computation of the transfer characteristics of the fractional-order
dynamical systems has been the subject of several publications, e.g.
by numerical methods \cite{Dorcak2}, as well as analytical
methods \cite{Podlubny1}.
In the synthesis of a controller its parameters are determined according
to the given requirements. These requirements are, for example, the stability
measure, the accuracy of the regulation process, dynamical properties etc.
The check whether the requirements have been met can be done with a simulation
on the control circuit model. There are a large number of methods for  the
design of integer-order controllers, but the situation is worse in the case
of fractional-order controllers where the methods are only being worked out.
One of the methods being developed is the method (modification
of roots locus method) of dominant roots \cite{Petras1}, based on the given
stability measure and the damping measure of the control circuit.

\section{Basic mathematical tools for fractional calculus}

The fractional calculus is a generalization of integration and derivation
to non-integer order operators. The idea of fractional calculus has been
known since the development of the regular calculus, with the first reference
probably being associated with Leibniz and L'Hospital in 1695.
At first, we generalize the differential and integral operators into one
fundamental operator $D^{\alpha}_{t}$  which is known as fractional
calculus:
$$
  _{a}D^{\alpha}_{t} = \left \{
        \begin{array}{ll}
                \frac{d^{\alpha}}{dt^{\alpha}} & \mbox{$\Re(\alpha)>0$,} \\
                 1 & \mbox{$\Re(\alpha)=0$,} \\
                \int_{a}^{t} (d\tau)^{-\alpha} & \mbox{$\Re(\alpha)<0$.}
        \end{array}
        \right.
$$
The two definitions used for the general fractional differintegral
are the Gr\"unwald  definition and the Riemann-Liouville (RL)
definition \cite{Oldham}. The Gr\"unwald definition is given here
\begin{equation}\label{GD}
     _{a}D^{\alpha}_{t}f(t)=\lim_{h \to 0}\frac{1}{h^{\alpha}}
                            \sum_{j=0}^{[\frac{t-a}{h}]}(-1)^j
                                       {\alpha \choose j} f(t-jh),
\end{equation}
where $[x]$ means the integer part of $x$. The RL definition is given as
\begin{equation}\label{LRL}
   _{a}D_{t}^{\alpha}f(t)=
    \frac{1}{\Gamma (n -\alpha)}
    \frac{d^{n}}{dt^{n}}
    \int_{a}^{t}
    \frac{f(\tau)}{(t-\tau)^{\alpha - n + 1}}d\tau, \\
\end{equation}
for $(n-1 < \alpha <n)$ and
where $\Gamma (x)$ is the well known Euler's \textit{gamma} function.

The Laplace transform method is used for solving engineering problems.
The formula for the Laplace transform of the RL fractional derivative
(\ref{LRL}) has the form \cite{Podlubny1}:
\begin{eqnarray}
      \int_{0}^{\infty} e^{-pt}\, _{0}D_{t}^{\alpha}f(t) \, dt =
      \nonumber \\
       = p^{\alpha}F(p) - \sum_{k=0}^{n-1} p^{k} \,
         \left.  _{0}D_{t}^{\alpha-k-1}f(t) \right|_{t=0},
\end{eqnarray}
for $(n-1 < \alpha \leq n)$.

For numerical calculation of fractional-order derivation we can use the
relation (\ref{FD}) derived from the Gr\"unwald definition (\ref{GD}).
This relation has the following form:
\begin{equation}\label{FD}
    _{(t-L)}D^{\alpha}_tf(t) \approx
      h^{-\alpha}  \sum_{j=0}^{N(t)} b_jf(t-jh),
\end{equation}
where $L$ is the "memory length", $h$ is the step size of the
calculation,
$$
   N(t) = \mbox{min}\left\{\left [\frac{t}{h}\right ],
             \left [\frac{L}{h}\right ]\right\},
$$
[x] is the integer part of $x$ and $b_j$ is the binomial
coefficient:
\begin{equation}\label{b_k}
  b_0 = 1, \qquad b_j = \left (1 -
                              \frac{1+\alpha}{j}\right)b_{j-1}.
\end{equation}

For the solution of the fractional-order
differential equations (FODE)  most effective and easy
analytic methods were developed based on the formula of the Laplace
transform method of the Mittag-Leffler function in two parameters
\cite{Podlubny1}.
A~two-parameter function of the Mittag-Leffler type is defined
by the series expansion:
\begin{equation} \label{ML-Definition}
E_{\alpha, \beta}(z) = \sum_{k=0}^{\infty}\frac{z^k}{\Gamma (\alpha k+\beta)},
\hspace{1em} (\alpha, \beta > 0).
\end{equation}
The Mittag-Leffler function is a generalization of exponential
function $e^{z}$ and the exponential function is a particular
case of the Mittag-Leffler function. Here is the relationship given in
\cite{Podlubny1}:
$$
E_{1,1}(z) = \sum_{k=0}^{\infty}\frac{z^k}{\Gamma (k + 1)} =
             \sum_{k=0}^{\infty}\frac{z^k}{k!}=e^{z}.
$$

\section{Fractional - order closed \protect \\ control loop}

We will be studying feed-back control system with unit gain in the
feed-back loop (Fig.1), where $G_{r}(p)$ is the controller transfer
function, $G_{s}(p)$ is controlled system transfer function,
$W(p)$ is an input, $E(p)$ is an error, $U(p)$ is output from controller
and $Y(p)$ is output from system.
\begin{figure}[h]
           \vskip 3 mm
	   \centerline{\includegraphics[width=7.9cm]{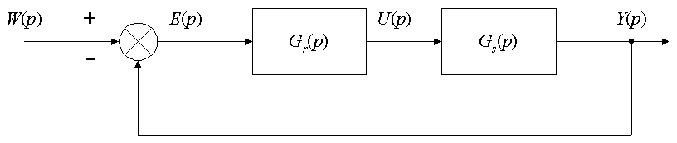}}
           \vskip 1 mm
           \centerline{Fig.1: {\it Feed - back control circuit}}
           \vskip 2 mm
\end{figure}
The transfer function of close feed-back control circuit (Fig.1)
has the form:
\begin{equation}\label{G_c}
        G_c(p)=\frac{Y(p)}{W(p)}=
               \frac{G_r(p)G_s(p)}{1+G_r(p)G_s(p)}.
\end{equation}

\subsection{Fractional - order controlled \protect \\ systems}

The fractional-order controlled system will be represented with fractional
model with transfer function  given by the following expression
\cite{Podlubny2}:
\begin{equation}\label{Gs}
      G_{s}(p) =\frac{Y(p)}{U(p)}=
                \frac{1}{a_{n}p^{\beta_{n}} + \ldots
                 + a_{1}p^{\beta_{1}} + a_{0}p^{\beta_{0}}},
\end{equation}
where  $\beta_k$, $(k = 0, 1, 2, \dots, n)$ are generally real numbers,
$\beta_{n} > \beta_{n-1} > \ldots > \beta_{1} > \beta_{0} \geq 0$
and $a_k$ $(k= 0, 1, \dots, n)$ are arbitrary constants.
In the  time  domain, the  transfer  function (\ref{Gs}) corresponds to
the $n$\/-term FODE with constant coefficients
\begin{equation} \label{n-term-equation}
     a_{n}\, D^{\beta_{n}}_ty(t) + \ldots +
     a_{1}\, D^{\beta_{1}}_ty(t) +
     a_{0}\, D^{\beta_{0}}_ty(t) = u(t).
\end{equation}

Identification   methods \cite{Dorcak1,Podlubny3} for
determination of the coefficients $a_k$ and  $\beta_k$, $(k= 0, 1, \dots, n)$
were developed, based on minimization of the difference between
the calculated ($y^{c}$) and experimentally measured ($y^{e}$) values
$$
        Q=\frac{1}{M+1}\sum_{m=0}^{M}[y^{e}_{m}-y^{c}_{m}]^2,
$$
where $M$ is the number of measured values.

For the analytical solution of the  $n$\/-term FODE
(\ref{n-term-equation})
we can write formula in general form \cite{Podlubny1}:
\begin{eqnarray}
     y(t)  =
     \nonumber \\
     =  \frac{1}{a_n} \sum_{m=0}^{\infty} \frac{(-1)^m}{m!}
     \sum_{{k_0+k_1+\dots +k_{n-2}=m
                  \atop
             k_0 \geq 0; \dots, k_{n-2} \geq 0
            }}
     (m; k_0, k_1, \dots , k_{n-2})
     \nonumber \\
    \prod_{i=0}^{n-2}
     \left(
          \frac{a_i}{a_n}
     \right)^{k_i}
     t^{(\beta_n -\beta_{n-1})m +\beta_n
                      +\sum_{j=0}^{n-2}(\beta_{n-1}-\beta_j)k_j -1}
     \nonumber \\
     E^{(m)}_{\beta_n-\beta_{n-1}, +\beta_n
                      +\sum_{j=0}^{n-2}(\beta_{n-1}-\beta_j)k_j }
         \left(
               -\frac{a_{n-1}}{a_n}t^{\beta_n-\beta_{n-1}}
         \right),
      \nonumber
\end{eqnarray}
where $E_{\lambda, \mu}(z)$ is the
Mittag-Leffler function in two parameters (\ref{ML-Definition}),
$$
E_{\lambda,\mu}^{(n)}(y) \equiv \frac{d^{n}}{dy^{n}}E_{\lambda ,\mu}(y) =
\sum_{j=0}^{\infty} \frac{(j+n)! \,\, y^{j}}
                         {j! \,\, \Gamma (\lambda j + \lambda n + \mu)},
$$
for $(n = 0, 1, 2, ...)$.

For the numerical solution of the $n$\/-term FODE (\ref{n-term-equation})
we can write formula in general form:
$$
y(k)=\frac{u(k)-\sum_{i=1}^{n}(a_ih^{-\beta_i}\sum_{j=1}^{k}b_jy(k-j))}
         {\sum_{i=0}^{n}a_ih^{-\beta_i}b_0},
$$
for $k=1,2,3,\dots$, $y(0)=0$, $y(1)=0$, where $u(k)$ is a function on
the right side of the differential equation.

\subsection{Fractional - order controllers}

The fractional-order controller will be  represented  by
fractional-order $PI^{\lambda}D^{\delta}$ controller with transfer
function  given by the following expression \cite{Podlubny2}:
\begin{equation}\label{Gr}
      G_{r}(p)=\frac{U(p)}{E(p)}=K+T_ip^{-\lambda}+T_dp^{\delta},
\end{equation}
where $\lambda$ and $\delta$ are an arbitrary real numbers
$(\lambda, \delta \geq 0)$, $K$ is amplification (gain), $T_i$ is
integration constant and $T_d$ is differentiation constant. In the time
domain equation (\ref{Gr}) has the form:
\begin{equation}\label{controller}
    u(t)=Ke(t)+T_iD_t^{-\lambda}e(t)+T_dD_t^{\delta}e(t).
\end{equation}

Taking $\lambda=1$ and $\delta=1$, we obtain a classical $PID$
controller. If $\lambda = 0$ $(T_i = 0)$ we obtain a $PD^{\delta}$
controller, etc.
All these types of controllers are particular cases of the
$PI^{\lambda}D^{\delta}$ controller. The $PI^{\lambda}D^{\delta}$
controller (\ref{controller}) is more flexible and gives an
opportunity to better adjust the dynamical properties of
the fractional-order control system.

\section{Synthesis of fractional - \protect \\ order controllers}

For the design of fractional-order controller  a new
method was suggested based on the dominant roots principle.
This method is based on from poles distribution of the characteristic
equation in the complex plane (Fig.2).
\begin{figure}
    \vskip 2 mm
    \includegraphics[height=6cm]{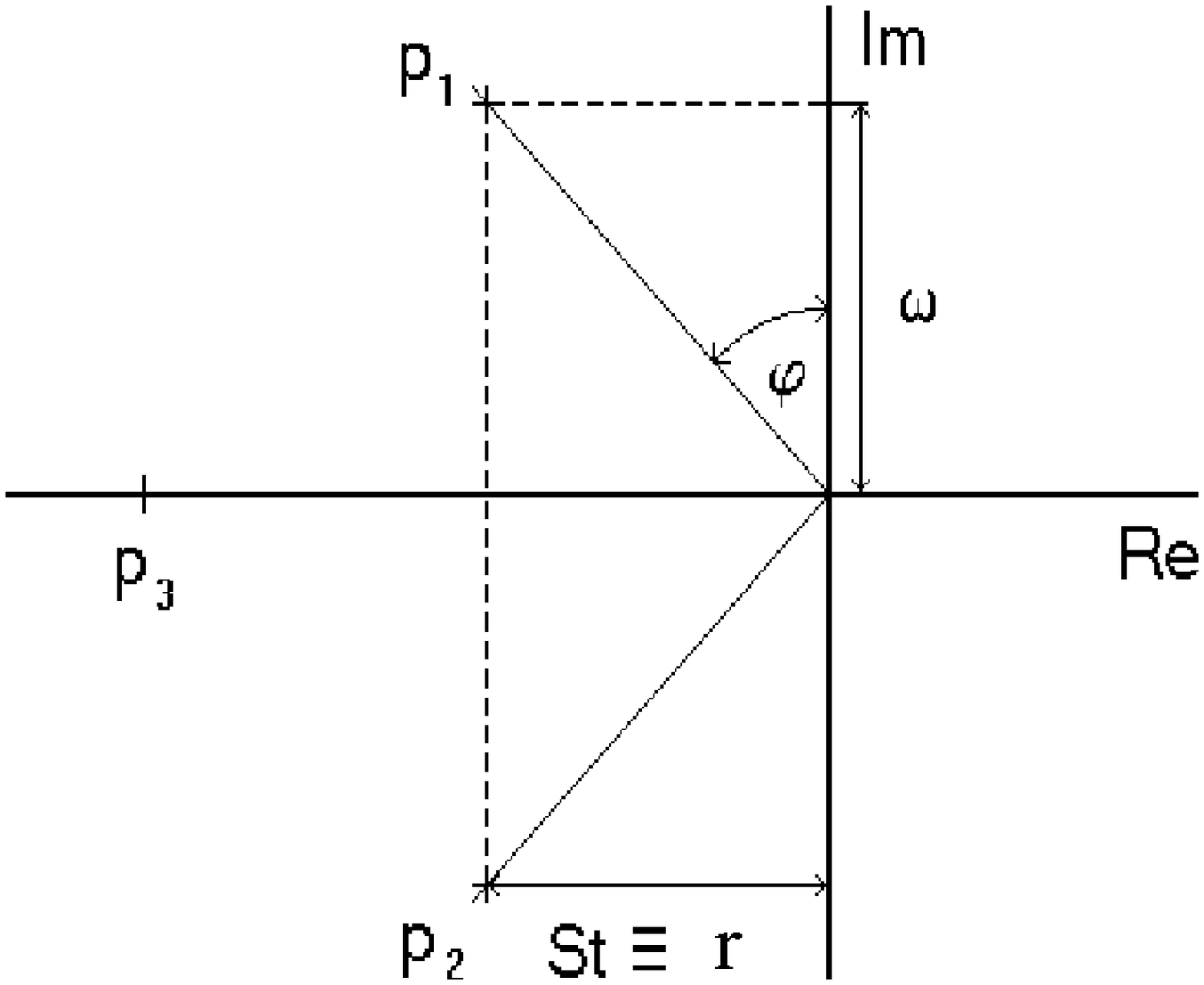}
    \centerline{Fig.2: {\it Roots in the complex plane}}
    \vskip 3 mm
\end{figure}

Values of dominant roots are designed for the quality requirement of the
control circuit. Their significance is, that dominant roots are defined for
stability measure $S_t$ and damping measure $T_l$. Under
there conditions the  complex conjugate roots satisfy the
equation:
\begin{equation}\label{komplex}
   p_{1,2} = -r \pm i \omega.
\end{equation}

The parameters of controller were set up so that other poles were from
dominant remote to the left side.
The parameters design of the fractional-order controller can be
divided into two stages:
\begin{enumerate}
\item \textbf{Design of parameter $K$} \\

Proportional parameter $K$ influence the value of static deviation
$E_t\,[\%]$, control time $T_r\,[s]$, and overshoot $P_r\,[\%]$.
Generally, with increased value of $K$, control time $T_r\,[s]$ is
decreasing and static deviation $E_t\,[\%]$ is lower:
$$
   K~\geq ({100}/{E_t}) - a_0.
$$

\item  \textbf{Design of parameters $T_d$, $\delta$, $T_i$, $\lambda$} \\

We define the required stability measure $S_t=r$ and damping
measure $T_l=\frac{r}{\omega}$. This requirement is satisfied
by the  complex conjugate roots (\ref{komplex}). We will similarly
use characteristic equation as the classical root locus method \cite{Dorf}.
The characteristic equation of fractional-order control loop
(\ref{G_c}) has the form:
\begin{equation}\label{ch_r}
  G_r(p)G_s(p) + 1 = 0.
\end{equation}

After substitution of the fractional-order controller transfer function
(\ref{Gr}) and the fractional-order controlled system transfer function
(\ref{Gs}) and after some manipulation we obtain characteristic equation
(\ref{ch_r}) in the following form:
\begin{equation}\label{ch_r_u}
\sum_{k=0}^{n}a_{k}p^{\beta_{k}}+(K+T_ip^{-\lambda}+T_dp^{\delta})=0.
\end{equation}

This algebraic equation is in the general form for the $n$\/-term
fractional-order controlled system and fractional-order
$PI^{\lambda}D^{\delta}$ controller. Solution of this equation for
known poles and parameter $K$ gives unknown parameters
$T_d$, $\delta$, $T_i$,
$\lambda$. The algebraic equation (\ref{ch_r_u}) is solved in
the space of complex variable. Result are parameters of
fractional-order controller for required stability measure and
damping measure.
\end{enumerate}

\section{Control algorithm for \protect \\ fractional - order  controllers}

The control algorithm was designed according to the control scheme
in Fig.1. The position algorithm uses discrete time steps \it
k~\rm and consists of the following steps \cite{Petras2}:
\begin{enumerate}
\item  \textbf {Filtering of required value:}  \\
       $$
          w^{*}(k) = w^{*}(k-1) + 0.5 (w(k) - w^{*}(k-1)),
       $$
       where $w(k)$ is required value.
\item  \textbf {Calculating the control error:} \\
       $$
          e(k)=w^{*}(k) - y(k),
       $$
       for discrete time step $(k=1,2,\dots)$,
       where $e(k)$ is regulation error and $y(k)$ is measured value.
\item  \textbf {Determination of control value:} \\
       \begin{eqnarray}
        u(k)= K~e(k) + \frac{T_i}{T^{-\lambda}}\sum\limits_{j=v}^{k}{q_je(k-j)} +
             \nonumber \\
              + \frac{T_d}{{T}^{\delta}}\sum\limits_{j=v}^{k}{d_je(k-j)},
             \nonumber
        \end{eqnarray}
      for discrete time step $(k=1,2,\dots)$, where $T$ is the length
      of time step (sample period).
      The binomial coefficients $d_j$ and $q_j$ were calculated
      from the recurrent equation (\ref{b_k}).
      The numerical algorithm (\ref{FD})
      requires to store the whole history
      $(v=0)$.
      For improving their effectiveness we have used "short memory"
      principle \cite{Dorcak2}, where $v=0$ for $k<(L/T)$ or $v=k-(L/T)$
      for $k>(L/T)$.
      Besides the "short memory" the control quality is influenced
      by time step $T$.
\end{enumerate}
Fractional-order controllers can be realized as software or by passive
or active electrical elements \cite{Ichise, Nakagawa}.

\section{Example}

We verify the above methods on an example from \cite{Dorcak2}.
Assume the system which we can describe by three-term $(n=2)$
differential equation  with coefficients $a_2 = 0.8, a_1 = 0.5,
a_0 = 1, \beta_2 = 2.2, \beta_1 = 0.9, \beta_0=0$. After its approximation
with an integer-order system we have $a_2=0.7414, a_1=0.2313, a_0=1,
\beta_2 = 2, \beta_1 = 1, \beta_0=0$. The integer-order  $PD$
controller  designed with the dominant roots method applied to the
approximated system has the parameters $K=20.5, T_d=2.7343, \delta =1$.
By applying the controller to the original system we do not achieve
the required quality of the regulation process as with the approximated
integer-order system. This was confirmed via simulation in the time
domain \cite{Dorcak2} and also by checking the stability measure
and damping measure.

This proves the inadequacy of approximating the fractional-order
system with an integer-order system for the purpose of controller design.
It is suitable to consider fractional-order systems and also
controllers, by means of which it is possible to obtain higher quality
regulation and robustness. For the fractional-order
$PD^{\delta}$ controller we then have $K=20.5, T_d=5.79, \delta=0.95$ and
the required quality of regulation is ensured also in the original
fractional-order system \cite{Petras1}.
\begin{figure}[h]
           \vskip 3 mm
           \centerline{\includegraphics[width=8cm]{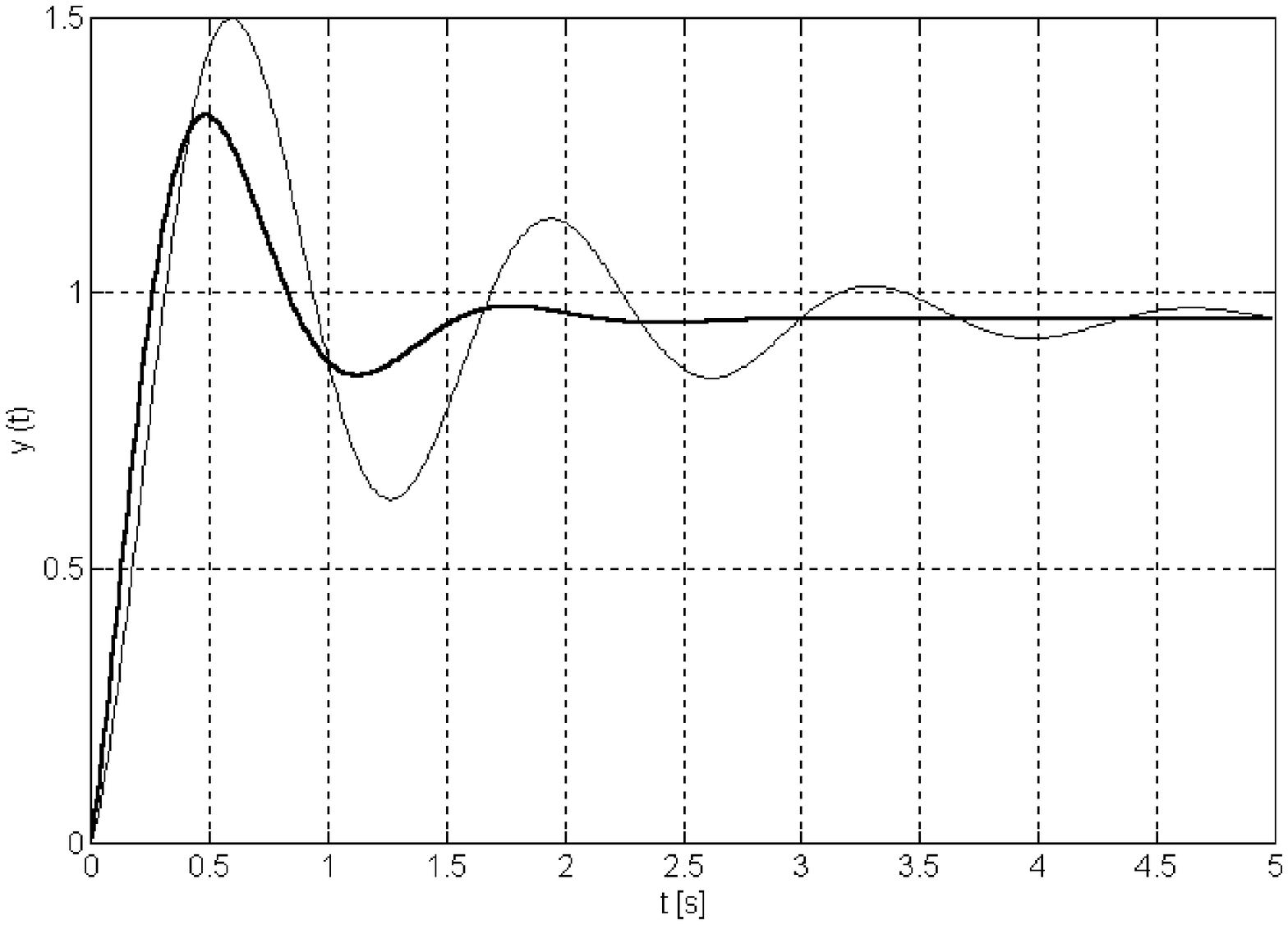}}
           \centerline{Fig.3: {\it Transient unit - step responses}}
           \vskip 2 mm
\end{figure}

On the Fig.3 is showed the comparison of the unit-step responses
the integer-order $PD$ controller applied on the
fractional-order system (thin line) and the fractional-order
$PD^{\delta}$ controller applied on the fractional-order
system (thick line).

The dynamical properties of the closed loop with fractional-order
controlled system and the fractional-order controller are better
than the dynamical properties of the closed loop with the
integer-order controller. The systems with the integer-order
controller stabilizes slower and has larger surplus oscillations.
We can see that use of the fractional-order controller leads to
the improvement of the control of the fractional-order system.
\section{Conclusion}

The above methods make it possible to design fractional-order
controllers with given measures of stability and damping.
The results of previous works also show that fractional-order
controllers are more robust \cite{Petras3}, which means they are
less sensitive to changes of the system parameters and controller
parameters.
This can even lead to qualitatively different dynamical phenomena
in control circuits.

\section*{Acknowledgement}

The author thanks to Dr. Lubomir Dorcak  for discussion of this
problem and also to Dr. Ladislav Pivka for checking the
language of this paper.
This work was supported by grant VEGA 1/4333/97 from
the Slovak Agency for Science.

\vskip 3mm
\textbf{Ivo Petras} was born in Kosice, Slovak Republic, in 1973.
He received the MSc. and Ph.D degrees in process control from
the Department of Informatics \& Process Control, Faculty B.E.R.G.
of the Technical University of Kosice. His main research interests
include fractional calculus in robust control and process control.
\end{document}